\documentclass{amsart}

\usepackage[all]{xypic}
\usepackage[centertags]{amsmath}
\usepackage{amsfonts}
\usepackage{epsf}
\usepackage{amscd}
\usepackage{amssymb}
\usepackage{amsthm}
\usepackage{newlfont}
\usepackage{amsxtra}
 \newtheorem{thm}{Theorem}[section]
 \newtheorem{cor}[thm]{Corollary}
 \newtheorem{lem}[thm]{Lemma}
 \newtheorem{prop}[thm]{Proposition}
 \theoremstyle{definition}
 \newtheorem{defn}[thm]{Definition}
 \theoremstyle{remark}
 \newtheorem{rem}[thm]{Remark}
 \theoremstyle{definition}
 \newtheorem{conj}[thm]{Conjecture}
 \theoremstyle{remark}
 \newtheorem{example}[thm]{Example}
 \theoremstyle{definition}
 \newtheorem{notn}[thm]{Notation}
 \numberwithin{equation}{section}
 \newcommand{\SL}{\mathrm{SL}}

 \newcommand{\Type}{\mathrm{Type}}

 \newcommand{\Ver}{\mathrm{Ver}}
 \newcommand{\St}{\mathrm{St}}
 \newcommand{\Lk}{\mathrm{Lk}}
 
 \newcommand{\fm}{\mathfrak m}

 \newcommand{\fT}{\mathfrak T}

 \newcommand{\fB}{\mathfrak B}
 \newcommand{\cK}{\mathcal{K}}

 \newcommand{\cG}{\mathcal{G}}

 \newcommand{\cF}{\mathcal{F}}

 \newcommand{\R}{\mathbb{R}}
 
 \newcommand{\F}{\mathbb{F}}

 \newcommand{\eps}{\varepsilon}

 \newcommand{\G}{\Gamma}
 
 \newcommand{\La}{\Lambda}
 \newcommand{\la}{\lambda}

\begin{document}

\title[On the eigenvalues of $p$-adic curvature]
{On the eigenvalues of $p$-adic curvature}

\author{Mihran Papikian}

\address{Department of Mathematics, Pennsylvania State University, University Park, PA 16802}

\email{papikian@math.psu.edu}

\subjclass{Primary 20E42; Secondary 05C50,  22E40}

\keywords{Laplacians, buildings, discrete subgroups of $p$-adic
groups}

\date{}
%%% ----------------------------------------------------------------------

\begin{abstract}
We determine the maximal eigenvalue of the $p$-adic curvature
transformations on Bruhat-Tits buildings, and we give an essentially
optimal upper bound on the minimal non-zero eigenvalue of these
transformations.
\end{abstract}

% ----------------------------------------------------------------------
\maketitle
% ------------------------------------------------------------------------

\section{Statement of the results}

Let $\cK$ be a non-archimedean locally compact field with finite
residue field of order $q$. Let $G$ be an almost simple linear
algebraic group defined over $\cK$ of $\cK$-rank $\ell+1$. Let $\fT$
be the Bruhat-Tits building associated with $G(\cK)$ \cite{BT}. This
is an infinite, locally finite, contractible simplicial complex of
dimension $\ell+1$. Let $X$ be the link of a vertex of $\fT$. $X$ is
a finite simplicial complex of dimension $\ell$, which is a building
in the sense of Tits \cite{Brown}. In \cite{Garland}, Garland
defined a certain combinatorial Laplace operator $\Delta$ acting on
the $i$-cochains $C^i(X)$, $0\leq i\leq \ell-1$; see Definiton
\ref{defn1.7}. $\fT$ can be realized as the skeleton of a
non-archimedean symmetric space \cite[Ch. 5]{Berkovich}, and from
this point of view the operators $\Delta$ are the non-archimedean
analogues of curvature transformations of riemannian symmetric
spaces. Denote by $m^i(X)$ the minimal non-zero eigenvalue of
$\Delta$ acting on $C^i(X)$. By a rather ingenious argument, Garland
proved that for any $\eps>0$ there is a constant $q(\eps, \ell)$
depending only on $\eps$ and $\ell$ such that if $q>q(\eps, \ell)$
then $m^i(X)\geq \ell-i-\eps$. Denote by $M^i(X)$ the maximal
eigenvalue of $\Delta$. The main result of this paper is the
following (see Theorems \ref{prop_ny} and \ref{thm-last}):
\begin{thm}
$M^i(X)=\ell+1$ and $m^i(X)\leq \ell-i$.
\end{thm}
In fact we prove this result for an arbitrary finite building $X$.
Note that our estimate on $m^i(X)$ is the best possible estimate
which does not depend on $q$. Based on some explicit calculations,
we also propose a conjectural description of the behavior of all the
eigenvalues of $\Delta$ as $q\to \infty$; see Conjecture \ref{conj}.

The method of our proof is based on a modification of Garland's
original arguments. The results in \cite{Garland} are stated for
buildings. On the other hand, as is nicely explained in
\cite{Borel}, part of the argument in \cite{Garland} works for quite
general simplicial complexes. In  $\S$\ref{SecGM} we follow
\cite{Borel}.

The main application of Garland's estimate on $m^i(X)$ is a
vanishing result for the cohomology groups of discrete cocompact
subgroups of $G(\cK)$; see $\S$\ref{ssGV}. This vanishing theorem
plays an important role in many problems arising in representation
theory and arithmetic geometry. Incidentally, our explicit
calculations of the eigenvalues of Laplace operators indicate that,
despite the hope expressed in \cite{Garland}, Garland's method is
not powerful enough to prove the vanishing theorem unconditionally,
i.e., without a restriction on $q$ being sufficiently large.

\section{Proofs}\label{SecGM}

\subsection{Simplicial complexes} We start by fixing the terminology
and notation related to simplicial complexes.

A \textit{simplicial complex} is a collection $X$ of finite nonempty
sets, such that if $s$ is an element of $X$, so is every nonempty
subset of $s$. The element $s$ of $X$ is called a \textit{simplex}
of $X$; its \textit{dimension} is $|s|-1$. Each nonempty subset of
$s$ is called a \textit{face} of $s$. A simplex of dimension $i$
will usually be referred to as $i$-simplex. The \textit{dimension}
$\dim(X)$ of $X$ is the largest dimension of one of its simplices
(or is infinite if there is no such largest dimension). A
subcollection of $X$ that is itself a complex is called a
\textit{subcomplex} of $X$. The \textit{vertices} of the simplex $s$
are the one-point elements of the set $s$.

Let $s$ be a simplex of $X$. The \textit{star} of $s$ in $X$,
denoted $\St(s)$, is the subcomplex of $X$ consisting of the union
of all simplices of $X$ having $s$ as a face. The \textit{link} of
$s$, denoted $\Lk(s)$, is the subcomplex of $\St(s)$ consisting of
the simplices which are disjoint from $s$. If one thinks of $\St(s)$
as the ``unit ball'' around $s$ in $X$, then $\Lk(v)$ is the ``unit
sphere'' around $s$.

A specific ordering of the vertices of $s$ up to an even permutation
is called an \textit{orientation} of $s$. An \textit{oriented}
simplex is a simplex $s$ together with an orientation of $s$. Denote
the set of $i$-simplices by $\widehat{S}_i(X)$, and the set of
oriented $i$-simplices by $S_i(X)$. We will denote the vertices
$\widehat{S}_0(X)=S_0(X)$ of $X$ also by $\Ver(X)$. For $s\in
S_i(X)$, $\bar{s}\in S_i(X)$ denotes the same simplex but with
opposite orientation. An $\R$-valued \textit{$i$-cochain} on $X$ is
a function $f$ from the set of oriented $i$-simplices of $X$ to
$\R$, such that $f(s)=-f(\bar{s})$. Such functions are also called
\textit{alternating}. The $i$-cochains naturally form a $\R$-vector
space which is denoted $C^i(X)$. If $i<0$ or $i>\dim(X)$, we let
$C^i(X)=0$.

\subsection{Laplace operators} From now on we assume that $X$ is a finite $n$-dimensional
complex such that
\begin{enumerate}
\item[($1_X$)] Each simplex of $X$ is a face of some $n$-simplex.
\end{enumerate}
For $s\in S_i(X)$, let $w(s)$ be the number of (non-oriented)
$n$-simplices containing $s$. In view of ($1_X$), $w(s)\neq 0$ for
any $s$.

\begin{lem}\label{lem-w} Let $\sigma\in S_i(X)$ be fixed. Then
$$
\sum_{\substack{s\in \widehat{S}_{i+1}(X)\\ \sigma\subset
s}}w(s)=(n-i)\cdot w(\sigma).
$$
\end{lem}
\begin{proof}
Given an $n$-simplex $t$ such that $\sigma\subset t$ there are
exactly $(n-i)$ simplices $s$ of dimension $(i+1)$ such that
$\sigma\subset s\subset t$. Hence in the sum of the lemma we count
every $n$-simplex containing $\sigma$ exactly $(n-i)$ times.
\end{proof}

Define a positive-definite pairing on $C^i(X)$ by
\begin{equation}\label{eq-pairing}
(f,g):=\sum_{{s} \in \widehat{S}_i(X)} w(s)\cdot f(s)\cdot g(s),
\end{equation}
where $f,g\in C^i(X)$ and in $w(s)\cdot f(s)\cdot g(s)$ we choose
some orientation of $s$. (This is well-defined since both $f$ and
$g$ are alternating.)

Define the \textit{coboundary}, a linear transformation $d:
C^i(X)\to C^{i+1}(X)$, by
\begin{equation}\label{eq-d}
(df)([v_0,\dots,
v_{i+1}])=\sum_{j=0}^{i+1}(-1)^jf([v_0,\dots,\hat{v}_j,\dots,
v_{i+1}]),
\end{equation}
where $[v_0,\dots, v_{i+1}]\in S_{i+1}(X)$ and the symbol
$\hat{v}_j$ means that the vertex $v_j$ is to be deleted from the
array.

Let $s=[v_0,\dots, v_i]\in S_i(X)$ and $v\in \Ver(X)$. If the set
$\{v,v_0,\dots, v_i\}$ is an $(i+1)$-simplex of $X$, then we denote
by $[v, s]\in S_{i+1}(X)$ the oriented simplex $[v,v_0,\dots, v_i]$.
Define a linear transformation $\delta: C^i(X)\to C^{i-1}(X)$ by
\begin{equation}\label{eq-delta}
(\delta f)(s)=\sum_{\substack{v\in \Ver(X)\\ [v,s]\in S_i(X)}}
\frac{w([v,s])}{w(s)}f([v,s]).
\end{equation}
In (\ref{eq-d}) and (\ref{eq-delta}), by convention, an empty sum is
assumed to be $0$. One easily checks that $\delta$ is the adjoint of
$d$ with respect to (\ref{eq-pairing}):
\begin{lem}\label{v-prop1.12}
If $f\in C^i(X)$ and $g\in C^{i+1}(X)$, then $(df,g)=(f,\delta g)$.
\end{lem}

\begin{defn}\label{defn1.7} The \textit{Laplace operator} on $C^i(X)$ is the linear
operator $\Delta=\delta d$.
\end{defn}

Since $\Delta$ is self-adjoint with respect to the pairing
(\ref{eq-pairing}), and for any $f\in C^i(X)$, $(\Delta f,
f)=(df,df)\geq 0$, $\Delta$ is diagonalizable and its eigenvalues
are non-negative real numbers.

\begin{rem}
The Laplace operator in \cite[Def. 3.15]{Garland} is defined as
$\delta d+ d\delta$. What we denote by $\Delta$ in this paper is
denoted by $\Delta^+$ in \textit{loc. cit.} When $X$ is the link of
a vertex in a Bruhat-Tits building, Garland calls $\Delta^+$ the
\textit{p-adic curvature}; see \cite[p. 400]{Garland}.
\end{rem}

\subsection{Garland's method}\label{ssGM} For $v\in \Ver(X)$ let $\rho_v$ be
the linear transformation on $C^i(X)$ defined by:
$$
(\rho_vf)(s)=\left\{
  \begin{array}{ll}
    f(s) & \hbox{if } v\in s; \\
    0 & \hbox{otherwise.}
  \end{array}
\right.
$$
Since any $i$-simplex has $(i+1)$-vertices, for $f\in C^i(X)$ we
have the obvious equality
\begin{equation}\label{obv}
\sum_{v\in \Ver(X)}\rho_vf=(i+1)f.
\end{equation}
We also have the following obvious lemma:
\begin{lem}\label{v-prop1.13} \hfill
\begin{enumerate}
\item $\rho_v\rho_v=\rho_v$;
\item For $f \in C^i(X)$ and $g\in C^i(X)$, $(\rho_vf, g)=(f,
\rho_vg)$.
\end{enumerate}
\end{lem}

Let $d_v$ and $\delta_v$ be the linear operators $d$ and $\delta$
acting on the cochains of the finite simplicial complex $\Lk(v)$,
and let $\Delta_v:=\delta_vd_v$. Note that $\Lk(v)$ is an
$(n-1)$-dimensional complex satisfying condition $(1_X)$. For
$f,g\in C^i(X)$ define their inner product on $\Lk(v)$ by
\begin{equation}\label{eq-locp}
(f,g)_v:=\sum_{s\in \widehat{S}_i(\Lk(v))}w_v(s)\cdot f(s)\cdot
g(s),
\end{equation}
where $w_v(s)$ is the number of $(n-1)$-simplices in $\Lk(v)$
containing $s$. This is simply the pairing (\ref{eq-pairing}) of the
restrictions of $f$ and $g$ to $\Lk(v)$.

\begin{lem}\label{lem-borel} If $f\in C^i(X)$, then
$$
i\cdot (\Delta f, f)+(n-i)(f,f)=\sum_{v\in \Ver(X)}(\Delta \rho_v f,
\rho_v f).
$$
\end{lem}
\begin{proof} See \cite[Lem. 1.3]{Borel}. In the proof it is crucial that the
inner product $(\cdot, \cdot)$ on $C^i(X)$ is defined using the
weights $w(s)$.
\end{proof}

\begin{cor}\label{cor1.14} Let $f\in C^i(X)$.
If there is a positive real number $\La$ such that $$(\Delta \rho_v
f, \rho_v f)\leq \La\cdot (\rho_v f, f)$$ for all $v\in \Ver(X)$,
then
$$i\cdot (\Delta f, f)\leq \left(\La\cdot (i+1)-(n-i)\right)(f,f).$$
\end{cor}
\begin{proof}
This follows from Lemma \ref{v-prop1.13}, Lemma \ref{lem-borel} and
(\ref{obv}).
\end{proof}

From now on we assume that $i\geq 1$. Define a linear transformation
$\tau_v:C^{i}(X)\to C^{i-1}(X)$ by
$$
(\tau_vf)(s)=\left\{
  \begin{array}{ll}
    f([v,s]) & \hbox{if $s\in S_{i-1}(\Lk(v))$;}\\
    0 & \hbox{otherwise.}
  \end{array}
\right.
$$

\begin{lem}\label{prop7.12}
For $f, g\in C^i(X)$, we have
$(\tau_vf,\tau_vg)_v=(\rho_vf,\rho_vg)$.
\end{lem}
\begin{proof} We have
$$
(\tau_vf,\tau_vg)_v =\sum_{\sigma\in
\widehat{S}_{i-1}(\Lk(v))}w_v(\sigma)\cdot \tau_vf(\sigma)\cdot
\tau_vg(\sigma).
$$
It is easy to see that there is a one-to-one correspondence between
the $n$-simplices of $X$ containing $[v,\sigma]$ and the
$(n-1)$-simplices of $\Lk(v)$ containing $\sigma$, so
$w_v(\sigma)=w([v,\sigma])$. Hence the above sum can be rewritten as
$$
\sum_{s\in \widehat{S}_{i}(\St(v))}w(s)\cdot (\rho_vf)(s)\cdot
(\rho_vg)(s).
$$
Since $\rho_v f$ is zero away from $\St(v)$, the sum can be extended
to the whole $\widehat{S}_{i}(X)$, so the lemma follows.
\end{proof}

\begin{lem}\label{prop7.14} For $f\in C^i(X)$, we have
$(\Delta \rho_v f,\rho_v f)=(\Delta_v\tau_vf,\tau_vf)_v$.
\end{lem}
\begin{proof}
By Lemma \ref{v-prop1.13} and Lemma \ref{prop7.12},
$$
(\Delta \rho_v f,\rho_vf)=(\rho_v\Delta \rho_v
f,\rho_vf)=(\tau_v\Delta \rho_vf,\tau_vf)_v.
$$
Expanding $\tau_v\Delta\rho_v f(s)$ for $s\in S_{i-1}(\Lk(v))$, one
easily checks that $\tau_v\Delta\rho_v f=\Delta_v\tau_vf$, which
implies the claim.
\end{proof}

\begin{lem}\label{lem2.17}
If $c$ is an eigenvalue of $\Delta_v$ acting on $C^{i-1}(\Lk(v))$
for some $v\in \Ver(X)$, then $c$ is also an eigenvalue of $\Delta$
acting on $C^i(X)$.
\end{lem}
\begin{proof}
Let $f\in C^{i-1}(\Lk(v))$ be such that $\Delta_vf=c\cdot f$. Define
$g\in C^i(X)$ as follows. If $s\in S_i(X)$ does not contain $v$ then
$g(s)=0$. If $s=[v,\sigma]$ then $g(s)=f(\sigma)$. In particular,
$\tau_v g=f$. We know that $\tau_v \Delta \rho_v g=\Delta_v\tau_v
g$. Obviously $\rho_v g=g$, so $\tau_v \Delta g=\Delta_v f=c\cdot
f=\tau_v(c\cdot g)$. This implies that $\Delta g=c\cdot g$.
\end{proof}

\begin{notn} Given a finite simplicial complex $Y$ satisfying $(1_X)$,
let $M^i(Y)$ and $m^i(Y)$ be the maximal and minimal non-zero
eigenvalues of $\Delta$ acting on $C^i(Y)$, respectively. Denote
$$
\la^i_{\max}(Y):= \max_{\substack{v\in \Ver(Y)}}M^i(\Lk(v)) \quad
\text{and}\quad \la^i_{\min}(Y):= \min_{\substack{v\in
\Ver(Y)}}m^i(\Lk(v)).
$$
\end{notn}

\begin{cor}\label{cor2.19} $M^i(X)\geq \la^{i-1}_{\max}(X)$ and
$m^i(X)\leq \la^{i-1}_{\min}(X)$.
\end{cor}

\begin{prop}\label{cor7.15} For $f\in C^i(X)$, we have
$$
(\Delta \rho_vf,\rho_vf)\leq \la^{i-1}_{\max}(X)\cdot(\rho_vf,f).
$$
\end{prop}
\begin{proof}
By Lemma \ref{prop7.14}, $(\Delta
\rho_vf,\rho_vf)=(\Delta_v\tau_vf,\tau_vf)_v$. Let $\{e_1, \dots,
e_h\}$ be an orthogonal basis of $C^{i-1}(\Lk(v))$ with respect to
$(\cdot , \cdot)_v$ which consists of $\Delta_v$-eigenvectors. Write
$\tau_vf=\sum_j a_j e_j$. Then
$$ (\Delta_v \tau_v f, \tau_v f)_v \leq M^{i-1}(\Lk(v)) \sum_{j=1}^h
a_j^2 (e_j, e_j)_v \leq \la_{\max}^{i-1}(X)\cdot (\tau_v f, \tau_v
f)_v.$$
On the other hand, by Lemma \ref{v-prop1.13} and Lemma
\ref{prop7.12}, $ (\tau_vf,\tau_vf)_v =(\rho_vf,\rho_vf)= (\rho_vf,
f)$.
\end{proof}

Denote by $\tilde{H}^i(\Lk(v), \R)$ the $i$th reduced simplicial
cohomology group of $\Lk(v)$.

\begin{thm}[Fundamental Inequality]\label{thmFI}
For $1\leq i\leq n-1$, we have
$$
i\cdot M^i(X)\leq (i+1)\cdot \la^{i-1}_{\max}(X)-(n-i).
$$
If $\tilde{H}^{i-1}(\Lk(v), \R)=0$ for every $v\in \Ver(X)$, then
$$
i\cdot m^i(X)\geq (i+1)\cdot \la^{i-1}_{\min}(X)-(n-i).
$$
\end{thm}
\begin{proof} Let $f\in C^i(X)$ be such that $\Delta f=M^i(X)\cdot f$.
Proposition \ref{cor7.15} implies that the assumption of Corollary
\ref{cor1.14} is satisfied with $\La=\la^{i-1}_{\max}(X)$. This
proves the first part. The second part is Garland's original
fundamental estimate \cite[$\S$5]{Garland}. For a proof see Theorem
1.5 in \cite{Borel}.
\end{proof}

\begin{notn} For $m\geq 1$, let $I_m$ denote the $m\times m$ identity matrix and let $J_m$
denote the $m\times m$ matrix whose entries are all equal to $1$.
The minimal polynomial of $J_m$ is $x(x-m)$.
\end{notn}

\begin{example}\label{Ex-P} Let $X$ be an $n$-simplex. We claim that the
eigenvalues of $\Delta$ acting on $C^i(X)$ are $0$ and $(n+1)$ for
any $0\leq i \leq n-1$. It is easy to see that $0$ is an eigenvalue,
so we need to show that the only non-zero eigenvalue of $\Delta$ is
$(n+1)$, or equivalently, $m^i(X)=M^i(X)=n+1$. First, suppose $i=0$.
Since for any simplex of $X$ there is a unique $n$-simplex
containing it, one easily checks that $\Delta$ acts on $C^0(X)$ as
the matrix $(n+1)I_{n+1}-J_{n+1}$. The only eigenvalues of this
matrix are $0$ and $(n+1)$. Now let $i\geq 1$. The link of any
vertex is an $(n-1)$-simplex, so by induction
$\la_{\min}^{i-1}(X)=\la^{i-1}_{\max}(X)=n$. Since the reduced
cohomology groups of a simplex vanish, the Fundamental Inequality
implies
$$
i\cdot M^i(X)\leq (i+1)n-(n-i)=i(n+1)
$$
and
$$
i\cdot m^i(X)\geq (i+1)n-(n-i)=i(n+1).
$$
Hence $(n+1)\leq m^i(X)\leq M^i(X)\leq (n+1)$, which implies the
claim.
\end{example}

% ------------------------------------------------------------------------

\subsection{Buildings}\label{sec3}
Let $G$ be a group equipped with a Tits system $(G,B,N,S)$ of rank
$\ell+1$. To every Tits system, there is an associated simplicial
complex $\fB$ of dimension $\ell$, called the \textit{building} of
$(G,B,N,S)$. For the definitions and basic properties of buildings
we refer to Chapters IV and V in \cite{Brown}. The simplices of
$\fB$ are in one-to-one correspondence with proper parabolic
subgroups of $G$. Assume from now on that $G$ is finite. Then $\fB$
is a finite simplicial complex satisfying $(1_X)$. Given a simplex
$s$ of $\fB$, it is known that $\Lk(s)$ is again a building
corresponding to a Tits system of rank $\ell-\dim(s)$.

We would like to estimate $M^i(\fB)$ and $m^i(\fB)$ for $0\leq i\leq
\ell-1$. This will be done inductively, using induction on $i$ and
$\ell$. The base of induction is the following lemma:

\begin{lem}\label{lem-n2.2}
If $\ell=1$ then $M^0(\fB)=2$, and $m^0(\fB)\leq 1$.
\end{lem}
\begin{proof} When $\ell=1$,
the eigenvalues of $\Delta$ acting on $C^0(\fB)$ were calculated by
Feit and Higman in \cite{FH}. The claim follows from these
calculations. See also Proposition 7.10 in \cite{Garland} when $\fB$
is of Lie type.
\end{proof}

Let $K$ be the fundamental chamber of $\fB$, i.e., the
$\ell$-simplex of $\fB$ corresponding to the Borel subgroup $B$ of
the given Tits system. Every simplex $s$ of $\fB$ can be transformed
to a unique face $s'$ of $K$ under the action of $G$. Label the
vertices of $K$ by the elements of $I_\ell:=\{0,1,\dots,\ell\}$, and
define $\Type(s)$ to be the subset of $I_\ell$ corresponding to the
vertices of $s'$. $G$ naturally acts on $\fB$ and this action is
type-preserving and strongly transitive; see \cite[$\S$V.3]{Brown}.
From this perspective one can think of $K$ as the quotient $\fB/G$.

\begin{lem}\label{lem-extend}
$\ell+1$ is an eigenvalue of $\Delta$ acting on $C^i(\fB)$. In
particular, $M^i(\fB)\geq \ell+1$.
\end{lem}
\begin{proof} Given a function $f\in C^i(K)$, we can
lift it (uniquely) to a $G$-invariant function $\tilde{f}\in
C^i(\fB)$ defined by $\tilde{f}(\tilde{s}):=f(s)$, where $\tilde{s}$
is any preimage of $s$ in $\fB$. As is explained in
\cite[$\S$4.2]{Borel}, we have $\widetilde{\Delta f}=\Delta
\tilde{f}$. Hence the claim follows from Example \ref{Ex-P}.
\end{proof}

\begin{rem}\label{rem3.12}
Note that Lemma \ref{lem2.17}, coupled with Lemma \ref{lem-extend},
implies that the integers $\ell+1, \ell,\cdots, \ell-i+1$ are always
present among the eigenvalues of $\Delta$ acting on $C^i(\fB)$. On
the contrary, $\ell-i$ is not necessarily an eigenvalue.
\end{rem}

Let $f\in C^0(\fB)$, and let $R$ be a fixed constant. For each
$\alpha\in I_\ell$, define a function $f_\alpha$ on the vertices of
$\fB$ by $f_\alpha(v)=f(v)$ if $\Type(v)\neq \alpha$ and
$f_\alpha(v)=R\cdot f(v)$ if $\Type(v)=\alpha$. Also, for a fixed
$\alpha\in I_\ell$ define a linear transformation $\rho_\alpha$ on
$C^0(\fB)$ by
$$
\rho_\alpha=\sum_{\Type(v)=\alpha}\rho_v.
$$
For $f\in C^0(\fB)$ and any $\alpha$, we have
\begin{equation}\label{eq-jan1}
(\rho_\alpha df_\alpha, df_\alpha)=(df_\alpha,
df_\alpha)-((1-\rho_\alpha)df, df),
\end{equation}
and
\begin{equation}\label{eq-jan2}
(\Delta\rho_\alpha df_\alpha, \rho_\alpha
df_\alpha)=((1-\rho_\alpha)df,df)
\end{equation}
The equations (\ref{eq-jan1}) and (\ref{eq-jan2}) are the equations
(3) and (6) in \cite[$\S$4.5]{Borel}, respectively.

\begin{lem}\label{lemd31}
Let $f\in C^0(\fB)$ and suppose $\Delta f=c\cdot f$. Then
$$
\sum_{\alpha\in I_\ell}(\Delta f_\alpha,
f_\alpha)=\left[(\ell-c)(R-1)^2+c(R^2+\ell)\right]\cdot (f,f).
$$
\end{lem}
\begin{proof} Fix some type $\alpha$ and
let $g\in C^0(\fB)$ be a function such that $g(v)=0$ if
$\Type(v)\neq \alpha$. Then $(\Delta g, g)=\ell\cdot (g,g)$. Indeed,
$$
(\Delta g, g)=(dg,dg)=\sum_{[x,v]\in
\widehat{S}_1(\fB)}w([x,v])(g(v)-g(x))^2
$$
$$
=\sum_{\Type(v)=\alpha}g(v)^2\sum_{x\in
\Ver(\Lk(v))}w([x,v])=\ell\sum_{\Type(v)=\alpha}w(v)\cdot
g(v)^2=\ell\cdot (g,g).
$$
(The middle equality on the previous line follows from Lemma
\ref{lem-w}.) If we apply this to $g=f_\alpha-f$, then we get
\begin{equation}\label{eq-d31}
(\Delta f_\alpha, f_\alpha)=\ell\cdot (f_\alpha,
f_\alpha)-2(\ell-c)(f_\alpha, f)+(\ell-c)(f,f).
\end{equation}
We clearly have
$$
\sum_{\alpha\in I_\ell}f_\alpha= (\ell+R)\cdot f\quad \text{and}
\quad \sum_{\alpha\in I_\ell}(f_\alpha, f_\alpha)=(\ell+R^2)\cdot
(f,f).
$$
Summing (\ref{eq-d31}) over all types and using the previous two
equalities, we get the claim.
\end{proof}

\begin{thm}\label{prop_ny}
$M^i(\fB)=\ell+1$.
\end{thm}
\begin{proof}
By Lemma \ref{lem-extend}, it is enough to show that $M^i(\fB)\leq
\ell+1$. We start with $M^0(\fB)$. Let $f\in C^0(\fB)$. Since the
vertices of any simplex in $\fB$ have distinct types, one easily
checks that
$$
\sum_{\Type(v)=\alpha}(\Delta\rho_v df_\alpha,
\rho_vdf_\alpha)=(\Delta \rho_\alpha df_\alpha, \rho_\alpha
df_\alpha),
$$
so by Proposition \ref{cor7.15}
\begin{equation}\label{eq-ny}
(\Delta \rho_\alpha df_\alpha, \rho_\alpha df_\alpha)\leq
\la^0_{\max}(\fB)\cdot (\rho_\alpha df_\alpha, \rho_\alpha
df_\alpha).
\end{equation}
Since for any $v\in \Ver(\fB)$, $\Lk(v)$ is a building of dimension
$\ell-1$, the induction on $\ell$ gives $\la^0_{\max}(\fB)=\ell$.
Combining this with (\ref{eq-ny}), (\ref{eq-jan1}) and
(\ref{eq-jan2}), we get
\begin{equation}\label{eq-dec5}
(1+\ell)\cdot ((1-\rho_\alpha)df,df)\leq \ell\cdot (df_\alpha,
df_\alpha).
\end{equation}
Now assume $\Delta f=c\cdot f$. Note that
\begin{align}\label{eq-ny3}
\sum_{\alpha\in I_\ell}(1-\rho_\alpha)df&=(\ell+1)df-\sum_{v\in
\Ver(\fB)}\rho_v df\\ \nonumber &=(\ell+1)df-2df=(\ell-1)df,
\end{align}
so summing the inequalities (\ref{eq-dec5}) over all types and using
Lemma \ref{lemd31}, we get
\begin{equation}\label{eq-last}
(\ell+1)(\ell-1)c\cdot (f,f)\leq \ell\cdot
\left[(\ell-c)(R-1)^2+c(R^2+\ell)\right]\cdot (f,f).
\end{equation}
If we put $R=(\ell-c)/\ell$, then (\ref{eq-last}) forces $c\leq
\ell+1$. In particular, $M^0(\fB)\leq \ell+1$.

Now let $i\geq 1$. The induction on $i$ and $\ell$ implies that
$\la^{i-1}_{\max}(\fB)=\ell$. From the Fundamental Inequality
\ref{thmFI} we get
$$
i\cdot M^i(\fB)\leq (i+1)\cdot \ell-(\ell-i),
$$
which implies $M^i(\fB)\leq \ell+1$.
\end{proof}

\begin{thm}\label{thm-last}
$m^i(\fB)\leq \ell-i$.
\end{thm}
\begin{proof} We start with $i=0$.
Denote $c:=m^0(\fB)$ and let $f$ be a $\Delta$-eigenfunction with
eigenvalue $c$. First we claim that $c\neq \ell+1$. Indeed, $\Delta$
is a semi-simple operator and if $c=\ell+1$ then by Theorem
\ref{prop_ny} it has only two distinct eigenvalues, namely $0$ and
$\ell+1$. This implies that $\Delta^2=(\ell+1)\Delta$. But it is
easy to check that this equality is false. Next, $c\neq \ell+1$
implies $(\Delta f_\alpha, f_\alpha)\geq c\cdot (f_\alpha,
f_\alpha)$; see equation (1) in \cite[$\S$4.6]{Borel}. Summing over
all types,
$$
\sum_{\alpha\in I_\ell} (\Delta f_\alpha, f_\alpha)\geq
c(\ell+R^2)\cdot (f, f).
$$
Comparing this inequality with the expression in Lemma \ref{lemd31},
we conclude that $$(\ell-c)(R-1)^2\geq 0.$$ Since $R$ is arbitrary,
we must have $c\leq \ell$.

Now assume $i\geq 1$. By Corollary \ref{cor2.19} and induction on
$i$ and $\ell$, we have $m^i(\fB)\leq \la^{i-1}_{\min}(\fB)\leq
(\ell-1)-(i-1)=\ell-i$.
\end{proof}

\begin{thm}[Garland]\label{thm-G} Assume that $G$ is the group of $\F_q$-valued
points of a simple, simply connected Chevalley group. For any
$\eps>0$ there is a constant $q(\eps, \ell)$ depending only on
$\eps$ and $\ell$, such that if $q>q(\eps, \ell)$ then $m^i(\fB)\geq
\ell-i-\eps$.
\end{thm}
\begin{proof} For the proof see Sections 6, 7, 8 in \cite{Garland},
or Proposition 5.4 in \cite{Borel}.
\end{proof}

\section{Examples}\label{SecEg}

In this section we compute explicitly in some cases the eigenvalues
of $\Delta$ acting on $C^i(\fB)$. We concentrate on
$G=\mathrm{SL}_{\ell+2}(\F_q)$ for small $\ell$, with $B\subset G$
being the upper triangular group and $N$ being the monomial group,
cf. \cite[$\S$V.5]{Brown}. Denote the corresponding building by
$\fB_{\ell,q}$. The dimension of $\fB_{\ell,q}$ is $\ell$. Denote by
$\fm^i_\ell(q;x)$ the minimal polynomial of $\Delta$ acting on
$C^i(\fB_{\ell,q})$, $0\leq i\leq \ell-1$.

First, we recall an elementary description of $\fB_{\ell,q}$ which
is convenient for actual calculations. Let $V$ be a linear space
over $\F_q$ of dimension $\ell+2$. A \textit{flag} in $V$ is a
nested sequence $\cF: F_0\subset F_1\subset \cdots \subset F_i$ of
distinct linear subspaces $F_0,\dots, F_i$ of $V$ such that $F_0\neq
0$ and $F_i\neq V$. $\fB_{\ell,q}$ is isomorphic to the simplicial
complex whose vertices correspond to the non-zero linear subspaces
of $V$ distinct from $V$; the vertices $v_0,\dots, v_i$ form an
$i$-simplex if the corresponding subspaces form a flag.

\vspace{0.1in}

Now assume $\ell=1$. In this case $\fB_{\ell,q}$ is isomorphic to
the $1$-dimensional complex whose vertices correspond to $1$ and
$2$-dimensional subspaces of a $3$-dimensional vector space $V$ over
$\F_q$, two vertices being adjacent if one of the corresponding
subspaces is contained in the other. With a slight abuse of
terminology, we will call $1$ and $2$ dimensional subspaces lines
and planes, respectively. The number of lines and planes in $V$ is
$m=q^2+q+1$ each. Let $A=(a_{ij})$ be the $m\times m$ matrix whose
rows are enumerated by the lines in $V$ and columns by the planes,
and $a_{ij}=-1$ if the $i$th line lies in the $j$th plane, and is
$0$ otherwise. We can choose a basis of $C^0(\fB_{\ell,q})$ so that
$(q+1)\Delta$ acts as the matrix
$$
(q+1)I_{2m}+\begin{pmatrix} 0 & A \\ A^t & 0\end{pmatrix}.
$$
Let $M=\begin{pmatrix} 0 & A \\ A^t & 0\end{pmatrix}$. Since any two
distinct lines lie in a unique plane and any line lies in $(q+1)$
planes, $AA^t=qI_m+J_m$. By a similar argument, $A^tA=qI_m+J_m$.
Hence
$$
M^2= qI_{2m}+\begin{pmatrix} J_m & 0 \\ 0 & J_m\end{pmatrix}.
$$
This implies that $(M^2-qI_{2m})(M^2-(q+1)^2I_{2m})=0$. Since
$(q+1)\Delta - (q+1)I_{2m}=M$, we conclude that $(q+1)\Delta$
satisfies the polynomial equation
$$
x(x-(2q+2))(x^2-(2q+2)x+(q^2+q+1))=0.
$$
It is not hard to see that this is in fact the minimal polynomial of
$(q+1)\Delta$. Hence
$$
\fm^0_1(q;x)=x(x-2)\left(x^2-2x+\frac{q^2+q+1}{q^2+2q+1}\right).
$$
The minimal non-zero root is $1-\sqrt{q}/(q+1)$. The smallest
possible value of this expression is approximately $0.53$, which
occurs at $q=2$, the value tends to $1$ as $q\to \infty$.

The very next case $\ell=2$ is already considerably harder to
compute by hand. With the help of a computer, we deduced that
\begin{align*}
\fm^0_2(q;x)=&x(x-2)\left(x-3\right)\left(x-\frac{2q^2+3q+2}{q^2+q+1}\right)\\
&\times\left(x^2-\frac{4q^2+3q+4}{q^2+q+1}x+\frac{4q^2+4}{q^2+q+1}\right).
\end{align*}
The minimal non-zero root is
$$
\frac{1}{2(q^2+q+1)}\left(4q^2+3q+4-\sqrt{8q^3+9q^2+8q}\right),
$$
which is at least $1.08$ and tends to $2$ from below as $q\to
\infty$. The whole polynomial tends coefficientwise to the
polynomial $x(x-3)(x-2)^4$ as $q\to \infty$. Next
\begin{align*}
\fm^1_2(q;x)=&x(x-1)(x-2)(x-3)\\
&\times\left(x^2-2x+\frac{q^2+1}{q^2+2q+1}\right)
\left(x^2-3x+\frac{2q^2+2q+2}{q^2+2q+1}\right)\\
&\times\left(x^2-4x+\frac{4q^2+6q+4}{q^2+2q+1}\right).
\end{align*}
The minimal non-zero root is $1-\sqrt{2q}/(q+1)$. This is always in
the interval $[1/3,1)$. Moreover, this eigenvalue is strictly larger
than $1/3$ for $q>2$ and tends to $1$ as $q\to \infty$; the whole
polynomial tends to $x(x-3)(x-1)^4(x-2)^4$.

The formulae for $\fm^0_2(q;x)$ and $\fm^1_2(q;x)$ are partly
conjectural, although almost certainly correct. We computed these
polynomials for $q=2, 3, 4, 5, 7$ using computer calculations with
concrete finite fields, and then came up with a formula which
recovers all the previous polynomials when we specialize $q$.

The complexity of calculations grows exponentially with $i$, $\ell$
and $q$, so for $\ell=3$ my computer was able to handle only $i=0$
for $q=2$ and $3$:
\begin{align*}
\fm^0_3(2;x)=& x(x - 4)\left(x - \frac{23}{7}\right)\left(x - \frac{19}{7}\right)\\
&\times \left(x^4 - 12x^3 + \frac{581528}{11025}x^2 -
\frac{220232}{2205}x + \frac{6734719}{99225}\right),
\end{align*}
\begin{align*}
\fm^0_3(3;x)=& x(x- 4)\left(x - \frac{42}{13}\right)\left(x - \frac{36}{13}\right)\\
&\times \left(x^4 - 12x^3 + \frac{14350977}{270400}x^2 -
\frac{2760633}{27040}x + \frac{309843369}{4326400}\right).
\end{align*}
The minimal non-zero roots of these polynomials are approximately
$1.68$ and $1.89$, respectively. To have a reasonable guess for the
coefficients of $\fm^0_3(q;x)$, one needs to compute these
polynomials for at least the next few values of $q$. Nevertheless,
note that the coefficients of above polynomials are close to the
coefficients of $x(x-4)(x-3)^6$.

The final example we have is
\begin{align*}
\fm^0_4(2;x)= & x(x - 4)(x - 5)\left(x - \frac{144}{35}\right)
\left(x^2 - \frac{1322}{155}x + \frac{2798}{155}\right)\\  &\times
\left(x^2 - \frac{276}{35}x + \frac{536}{35}\right) \left(x^3 -
\frac{1778}{155}x^2 + \frac{1306}{31}x - \frac{7512}{155}\right).
\end{align*}
The minimal non-zero root is approximately $2.32$, and the
coefficients of $\fm^0_4(2;x)$ are close to the coefficients of
$x(x-5)(x-4)^9$.

\vspace{0.1in}

The previous calculations, combined with Theorems
\ref{prop_ny}-\ref{thm-G} and Remark \ref{rem3.12}, suggest the
following possibility:

\begin{conj}\label{conj} In the situation of Theorem \ref{thm-G},
for any $\eps>0$ there is a constant $q(\eps, \ell)$ depending only
on $\ell$ and $\eps$ such that if $q>q(\eps, \ell)$ then any
non-zero eigenvalue of $\Delta$ acting on $C^i(\fB)$ is at a
distance less than $\eps$ from one of the integers $$\ell-i,\
\ell-i+1,\ \dots,\ \ell+1.$$
\end{conj}

\subsection{Garland's vanishing theorem}\label{ssGV}
Let $\cK$ be a field complete with respect to a non-trivial discrete
valuation and which is locally compact. Let $\F_q$ be the residue
field of $\cK$. Let $\cG$ be an almost simple linear algebraic group
over $\cK$. Suppose $\cG$ has $\cK$-rank $\ell+1$. Let $\fT$ be the
Bruhat-Tits building associated with $\cG(\cK)$. The link of a
simplex $s$ in $\fT$ is a finite building of dimension
$\ell-\dim(s)$. Using a discrete analogue of Hodge decomposition and
the Fundamental Inequality one proves the following theorem (see
\cite[Thm. 3.3]{Borel}):

\begin{thm}\label{thm4.2}
If $\la^{i-1}_{\min}(\fT)>\frac{\ell+1-i}{i+1}$, then $H^i(\G,
\R)=0$ for any discrete cocompact subgroup $\G$ of $\cG(\cK)$.
\end{thm}

Combining this with Theorem \ref{thm-G}, one concludes that there is
a constant $q(\ell)$ depending only on $\ell$ such that if
$q>q(\ell)$ then $H^i(\G, \R)=0$ for $1\leq i\leq \ell$. This is the
main result of \cite{Garland}. It is natural to ask whether the
restriction on $q$ being sufficiently large is redundant. This is
indeed the case, as was shown by Casselman \cite{Casselman}, who
proved the vanishing of the middle cohomology groups by an entirely
different argument.

Now let $\cG=\SL_{\ell+2}$. Then
$\la^{i-1}_{\min}(\fT)=m^{i-1}(\fB_{\ell, q})$ In all examples
discussed above $m^0(\fB_{\ell, q})>\ell/2$, so in these cases
Garland's method proves the vanishing of $H^1(\G,\R)$ without any
assumptions on $q$. On the other hand, $m^1(\fB_{2, 2})=1/3$. But to
apply Theorem \ref{thm4.2} to show that $H^2(\G,\R)=0$ we need
$\la_{\min}^1(\fT)>1/3$. Hence when $\ell=2$ we need to assume $q>2$
to conclude $H^2(\G,\R)=0$ from Garland's method.

% ------------------------------------------------------------------------

%\bibliographystyle{amsplain}
%\bibliography{pcurv}
\end{document}